\documentclass{amsart}
\usepackage{amssymb}
\newtheorem{thm}{Theorem}
\newtheorem{cor}{Corollary}
\newtheorem{lem}[thm]{Lemma}
\newtheorem{prop}[thm]{Proposition}
\newtheorem{defn}{Definition}

\newcommand{\set}[1]{\left\{#1\right\}}

\newcommand{\Complex}{\mathbb{C}}
\newcommand{\To}{\longrightarrow}
\newcommand{\image}[1]{\textrm{Im}(#1)}

\newcommand{\brck}[2]{\{#1,#2\}}
\newcommand{\jami}[2]{\sum\limits_{#1}^{#2}}
\newcommand{\all}[1]{\forall\,#1}
\newcommand{\vbar}{\Bigr|}
\newcommand{\imply}{\Rightarrow}
\newcommand{\one}{\mathbf{1}}
\newcommand{\smooth}[1]{C^\infty(#1)}
\newcommand{\fields}[1]{\mathfrak{F}(#1)}
\newcommand{\matr}[1]{M_{#1}(\Complex)}
\begin{document}
\title	[NC Symplectic Geometry of the Endomorphism Algebra]
	{Noncommutative Symplectic Geometry\\
	 of the Endomorphism Algebra\\
	 of a Vector Bundle}%
\author{Zakaria Giunashvili}%
\address{Department of Theoretical Physics,
         Institute of Mathematics,
         Georgian Academy of Sciences,
         Tbilisi, Georgia}%
\email{zaqro@gtu.edu.ge}%
\thanks{This work was supported by Nishina Memorial Foundation and was carried out during the author's visit in
Yokohama City University (Japan) in the framework of Nishina Memorial Foundation's Postdoctoral
Fellowship Program.}%
\subjclass{}%
\keywords{Noncommutative Geometry, Poisson Structure,
          Endomorphism Algebra, Noncommutative Submanifold,
          Noncommutative Quotient Manifold}%
\date{December 31, 2002}%
\begin{abstract}
We study noncommutative generalizations of such notions of the classical symplectic geometry as degenerate
Poisson structure, Poisson submanifold and quotient manifold, symplectic foliation and symplectic
leaf for associative Poisson algebras. We consider these structures for the case of the endomorphism algebra
of a vector bundle, and give the full description of the family of Poisson structures for this algebra.
\end{abstract}
\maketitle
\section{Introduction}
In this work we describe a noncommutative (algebraic) approach to such geometrical objects as
degenerate Poisson structures, Symplectic foliations for Poisson structures, symplectic
submanifolds and quotient manifolds. We generalize these notions for the case of noncommutative
differential calculus for an associative algebra. Our constructions are based on the
definitions of noncommutative submanifold and quotient manifold given in the work \cite{Masson1} of
Thierry Masson.

In \textbf{Section \ref{sectionA}} we give a brief overview of the definitions and some facts from the theory
of noncommutative submanifolds and quotient manifolds. Also, we introduce a notion of
\emph{locally proper} submanifold algebra which we need further as a noncommutative analogue of a
submanifold with $\textrm{codim}>0$.

\textbf{Section \ref{sectionB}} is devoted to the demonstration and investigation of noncommutative
submanifolds and quotient manifolds for the endomorphism algebra of a finite-dimensional complex
vector bundle. There are several works devoted to the noncommutative geometry of this algebra
(e.g., \cite{Masson2}, \cite{Dub-Masson}). Our purpose, in this section, is to give some
description of submanifold ideals for this algebra.

\textbf{Section \ref{sectionC}} is devoted to degenerate Poisson structures on general associative
algebra $A$. We introduce the definition of nondegenerate/degenerate Poisson structure and the
notion of locally proper Poisson ideal in $A$. Such ideal is a noncommutative analogue of
Poisson submanifold and when it is maximal it can be considered as a noncommutative symplectic leaf. The
existence of a locally proper Poisson ideal is some kind of indicator of
degeneracy of a Poisson structure. We study the relation between Poisson structures and
symplectic structures which is rather different from the classical case.

In \textbf{Section \ref{sectionD}} we apply the investigation and results of the previous sections
to the endomorphism algebra of a vector bundle: $End(E)$. The main result of this section is the
\textbf{full description of the family of Poisson structures on the endomorphism algebra of a vector
bundle, with $\dim(fiber)\geq 2$:} in this case \emph{any Poisson structure on $End(E)$ is of the form
$\brck{\Phi}{\Psi}=\lambda\cdot[\Phi,\Psi]$, where $\lambda$ is a smooth function on the base
manifold of the vector bundle}. The proof is mostly based on the fact that any Poisson structure on
the matrix algebra is of the form $\brck{a}{b}=k\cdot[a,b]$ for some constant number $k$ (the
proof of this statement is also given in this section), and the fact that the multiplicative commutator
$[End(E),End(E)]$ generates the entire algebra $End(E)$. It is worth to notice the general
proposition presented in this section which states that \textbf{for any such algebra $A$, the
commutator of which generates the entire algebra $A$ and any derivation vanishing on the center is
inner, a Hamiltonian map is always a $Z(A)$-module homomorphism from $A$ to $A/Z(A)$}, where
$Z(A)$ denotes the center of $A$. It seems that the center of an associative algebra plays very
important role for Poisson structures on this algebra and somehow ``fixes'' a Poisson structure.
In this section we give also the description of the symplectic leaves (in the noncommutative
sense) for the algebra $End(E)$.
\section[Noncommutative Submanifolds and Quotient Manifolds: a Brief Review]
{Submanifolds and Quotient Manifolds in Noncommutative Differential Calculus:
a Brief Review}\label{sectionA}
In this section we give a brief overview of definitions and facts from the theory of
noncommutative submanifolds and quotient manifolds in the derivation based noncommutative differential
calculus, which is based on the work of Thierry Masson (see \cite{Masson1}).

Let $A$ be an associative unital complex or real algebra. Let $I$ be a two-side ideal in
in $A$. We denote by $S$ the quotient algebra $A/I$ and by $p$ the quotient map from $A$ to $S$. We have
the following short exact sequence of algebra homomorphisms
\begin{equation}\label{shexseq_quotient-algebra}
0\To I\To A\stackrel{p}{\To}S\To 0
\end{equation}
Let $Der(A)$ be the Lie algebra of derivations of $A$. Introduce the following two Lie subalgebras
of $Der(A)$
$$
Der(A,I)=\set{X\in Der(A)\,\vbar\,X(I)\subset I}
$$
and its subalgebra
$$
Der(A,I)_0=\set{X\in Der(A)\,\vbar\,X(A)\subset I}
$$
The short exact sequence \ref{shexseq_quotient-algebra} induces the following exact sequence of Lie
algebra homomorphisms
$$
0\To Der(A,I)_0\To Der(A,I)\stackrel{\pi}{\To}Der(S)
$$
This exact sequence, in general, is not necessarily \textbf{short}.
\begin{defn}[see \cite{Masson1}] The quotient algebra $S=A/I$ is said to be a submanifold algebra
for $A$ if the homomorphism $\pi:Der(A,I)\To Der(S)$ is an epimorphism.
\end{defn}
In other words, in this case we have the short exact sequence
\begin{equation}\label{shexseq_vect-field-restriction}
0\To Der(A,I)_0\To Der(A,I)\stackrel{\pi}{\To}Der(S)\To 0
\end{equation}
The ideal $I$ is called the \emph{constraint ideal} for the submanifold algebra $S$. Sometimes, for brevity, we shall call $I$
a \emph{submanifold ideal}.

The classical geometric meaning of the Lie algebra $Der(A,I)$ is the following: it is the set of
such vector fields on a smooth manifold $M$ which are tangent to a submanifold $N\subset M$. It is
clear that if the submanifold $N$ is such that for each point $x\in N$ we have $T_x(N)=T_x(M)$,
then $Der(A,I)=Der(A)$. This is the case when $\textrm{codim}(N)=0$. For example, if $N$ is an open subset
in $M$. Translating this situation on the algebraic language, we introduce the following
\begin{defn}
We call a submanifold algebra \textbf{locally proper} if $Der(A,I)$ is a \textbf{proper} subalgebra
of $Der(A)$, i.e, $Der(A,I)\neq Der(A)$.
\end{defn}
Now, let us recall the noncommutative generalization of the notion of quotient manifold. Let $B$ be
a subalgebra of $A$. Define the following two Lie subalgebras of $Der(A)$
$$
T(A,B)=\set{X\in Der(A)\,\vbar\,X(B)\subset B}
$$
and its subalgebra
$$
V(A,B)=\set{X\in Der(A)\,\vbar\,X(B)=\{0\}}
$$
$V(A,B)$ is a Lie algebra ideal in $T(A,B)$. Moreover, $V(A,B)$ is the kernel of the Lie algebra
homomorphism $\rho:T(A,B)\To Der(B)$, which is just the restriction map.

Further, for any algebra $A$, we denote by $Z(A)$ the center of $A$.

\begin{defn}[see \cite{Masson1}] The subalgebra $B$ in $A$ is said to be a quotient manifold
algebra if the following three conditions are satisfied
\begin{itemize}
\item[q1. ] $Z(B)=B\cap Z(A)$
\item[q2. ] the restriction map $\rho: T(A,B)\To Der(B)$ is an epimorphism
\item[q3. ] $B=\set{b\in A\,\vbar\,X(b)=0,\,\all X\in V(A,B)}$
\end{itemize}
\end{defn}
\section[Submanifolds and Quotient Manifolds for the Endomorphism Algebra]
{Submanifolds and Quotient Manifolds for the Endomorphism Algebra of a Vector Bundle}
\label{sec_submanifolds-for-endomorphism-algebra}\label{sectionB}
In this section we investigate the endomorphism algebra of a finite-dimensional vector bundle in the
framework of the definitions and notions of the previous section.

Let $\pi:E\To M$ be a finite-dimensional complex (or real) vector bundle over a smooth manifold $M$.
We denote the algebra of endomorphisms of this bundle by $End(E)$. Any element $\Phi\in End(E)$ can
be considered as a section of the bundle of endomorphisms. For any point $x\in M$ we denote by $\Phi_x$
the value of this section at this point: $\Phi_x\in End(E_x)$, where $E_x=\pi^{-1}(x)$. If $\Phi$
is an element of the center of the algebra $End(E)$, then for any $x\in M$, the operator $\Phi_x$
is in the center of the algebra $End(E_x)$. But the algebra $End(E_x)$ is isomorphic to the algebra
of $n\times n$ matrices, where $n=\dim(E_x)$. The center of the matrix algebra is
$\Complex\cdot\mathbf{1}$, where $\mathbf{1}$ denotes the identity matrix. Hence we have that the
center of the algebra $End(E)$ is the set
$$
\set{f\cdot\mathbf{1}\,\vbar\,f\in C^\infty(M)}\equiv Z(End(E))
$$
For any associative algebra $A$ its center is invariant under the action of any derivation
of this algebra
$$
\begin{array}{c}
\all X\in Der(A),\ \all z\in Z(A),\ \all a\in A:\ za=az\imply X(za)=X(az)\imply\\
\\
\imply X(z)a+zX(a)=X(a)z+aX(z)\imply X(z)a=aX(z)\imply X(z)\in Z(A)
\end{array}
$$
Therefore we have the Lie algebra homomorphism
$$
\rho:Der(A)\To Der(Z(A))
$$
The Lie algebra of inner derivations of
$A$ is defined as
$$
Int(A)=\set{ad(a)=[a,\,\cdot]:A\To A\,\vbar\,a\in A}
$$
It is clear that $Int(A)\subset \ker{\rho}$.
Therefore, in the case when $\rho$ is an epimorphism and $\ker{\rho}=Int(A)$, we have that $Z(A)$
is a quotient manifold algebra for $A$ (see \cite{Masson1}).

Consider the above construction for the case when $A=End(E)$. In this case
$Der(Z(End(A)))=Der(C^\infty(M))=\mathfrak{F}(M)$, where $\mathfrak{F}(M)$ denotes the Lie algebra
of vector fields on the manifold $M$. Hence, we have a homomorphism
$$
\rho:Der(End(E))\To\mathfrak{F}(M)
$$
\begin{lem}\label{lem_internals-kernel}
The kernel of the homomorphism $\rho:Der(End(E))\To\mathfrak{F}(M)$ is the Lie algebra of inner
derivations of $End(E)$.
\end{lem}
\begin{proof}
Let $X\in\ker{\rho}$, then for any $f\in C^\infty(M)$ and $\Phi\in End(E)$, we have
$$
X(f\cdot\Phi)=f\cdot X(\Phi)
$$
which implies that $X:End(E)\To End(E)$ is a homomorphism of
$\smooth{M}$-modules. Therefore, $X$ corresponds to some homomorphism of the endomorphism bundle.
Hence, its action is pointwise
$$
X_m:End(E_m)\To End(E_m),\quad\all m\in M
$$
But as $End(E_m)$ is isomorphic to the matrix algebra, it has only inner derivations. And so $X$ is an inner
derivation.
\end{proof}

We denote by $\Gamma(E)$ the space of smooth sections of the vector bundle $(E,M,\pi)$. Let 
$$
\nabla:\Gamma(E)\To\Gamma(T^*(M))\otimes\Gamma(E)
$$
be a connection (covariant derivation)  on the
vector bundle $E$. For any $X\in\fields{M}$ define an operator $D_X:End(E)\To End(E)$ as:
\begin{equation}\label{form_connection-expand-endom}
D_X(\Phi)(s)=\nabla_X(\Phi(s))-\Phi(\nabla_X(s)),\quad\all\Phi\in End(E)\textrm{ and }\all s\in\Gamma(E)
\end{equation}
\begin{lem}
For any $X\in\fields{M}$ the operator $D_X$ is a derivation of the algebra $End(E)$ and
$\rho(D_X)=X$.
\end{lem}
\begin{proof}
For any pair of endomorphisms $\Phi,\Psi\in End(E)$ and $s\in\Gamma(E)$ we have the following
$$
\begin{array}{l}
D_X(\Phi\Psi)(s)=\nabla_X(\Phi(\Psi(s)))-(\Phi\Psi)(\nabla_X(s))=\\
\\
=D_X(\Phi)(\Psi(s))+\Phi(\nabla_X(\Psi(s)))-(\Phi\Psi)(\nabla_X(s))=\\
\\
=D_X(\Phi)(\Psi(s))+\Phi(D_X(\Psi)(s)+\Psi(\nabla_X(s)))-(\Phi\Psi)(\nabla_X(s))=\\
\\
=(D_X(\Phi)\circ\Psi+\Phi\circ D_X(\Psi))(s) 
\end{array}
$$
which implies that for any vector field $X\in\fields{M}$, the operator $D_X$ on $End(E)$ is a
derivation.

By definition of $\rho$, we have that for any $f\in\smooth{M}$
$$
(\rho(D_X))(f\cdot\mathbf{1})s=\nabla_X(fs)-f\nabla_X(s)=X(f)\cdot s
$$
which implies that $\rho(D_X)=X,\,\all X\in\fields{M}$.
\end{proof}

In fact, $D$ is the associated connection to $\nabla$ on the fiber bundle $End(E)$ (see \cite{Dub-Masson}).
It follows from the above lemma that $\rho$ is an epimorphism and this fact together with the
previous lemma implies that
\emph{the center of the endomorphism algebra $End(E)$, which is $\smooth{M}\cdot\mathbf{1}$, is a
quotient manifold algebra for $End(E)$}.

Hence, we have a short exact sequence of Lie algebra homomorphisms
\begin{equation}\label{shexseq_functions-derivs}
0\rightarrow Int(End(E))\rightarrow Der(End(E))\stackrel{\rho}{\rightarrow}Der(\smooth{M})\cong\fields{M}\rightarrow 0
\end{equation}
and the mapping $\fields{M}\ni X\mapsto D_X\in Der(End(E))$ gives a splitting of this short exact
sequence, but it is a splitting of the short exact sequence of $\smooth{M}$-module homomorphisms, because $X\mapsto D_X$
is a Lie algebra homomorphism only when the connection $\nabla$ is flat.

Further, in this section, we study the submanifold algebras for the algebra $End(E)$.

If $B\subset A$ is a quotient manifold algebra for an associative algebra $A$, and $I\subset B$ is a submanifold ideal in
$B$, then a natural candidate for a submanifold ideal in $A$ is the two-side ideal in $A$
generated by $I$. In general, such ideal is not always submanifold ideal, but this method gives a
positive result in some ``good'' cases. Let $N$ ($\partial N=\emptyset$) be a compact submanifold of $M$ and $I_N$ be the ideal
in $\smooth{M}$ consisting of the functions vanishing on $N$. The ideal in $End(E)$ generated by
$I_N$ coincides with the set of sections of the endomorphisms bundle vanishing on $N$. The quotient
algebra $End(E)/(I_N\cdot End(E))$ is canonically isomorphic to $End(E_N)$, where $E_N$ denotes the
restriction bundle of the vector bundle $E$ to the submanifold $N$:
$\pi_N=\pi|_N:E_N=\pi^{-1}(N)\To N$. The quotient map
$$
p:End(E)\To End(V_n)=End(E)/(I_N\cdot End(E))
$$
is just the restriction map, which maps any endomorphism of $\Phi:E\To E$ to its restriction
$\Phi_N:E_N\To E_N$.

Let $\nabla$ be a connection on $E$. If $X\in\fields{M}$ is such that $X(I_N)\subset I_N$, then
the ideal $I_N\cdot End(E)$ is invariant for the operator $D_X$, because:
$$
D_X(f\Phi)=X(f)\Phi+f\cdot D_X(\Phi)\in I_N\cdot End(E),\quad\all f\in I_N\textrm{ and }\all\Phi\in End(E)
$$
Therefore, the connection $\nabla$ defines a splitting of the exact sequence
$$
0\To Int(End(E_N))\To Der(End(E_N))\stackrel{\pi_N}{\To}\fields{N}\To 0
$$
This implies that any derivation $U_N\in Der(End(E_N))$ can be represented as
$$
U_N=ad(\Phi_N)+D_{X_N}
$$
where $\Phi_N\in End(E_N)$ and $X_N\in\fields{N}$. Let $X\in\fields{M}$
be an extension of $X_N$ and and $\Phi\in End(E)$ be an extension of $\Phi_N$, then we obtain that
the derivation $U=ad(\Phi)+D_X\in Der(End(E))$ is an extension of the derivation $U_N$. Hence the
derivations of $End(E_N)=End(E)/(I_N\cdot End(E))$ are ``covered'' by the derivations of $End(E)$,
which implies that $End(E_N)$ is a submanifold algebra for $End(E)$, and the ideal $I_N\cdot End(E)$
is the corresponding constraint ideal. Moreover, we have the following
\begin{prop}\label{prop_submanifold_in_base_total}
Let $I$ be an ideal in the algebra $End(E)$ and $I'$ be $I\cap(\smooth{M}\cdot\mathbf{1})$. If $I'$
is a submanifold ideal in the algebra $\smooth{M}\cdot\mathbf{1}$, corresponding to some compact closed submanifold $N$,
then $I$ is a submanifold ideal in the algebra $End(E)$
\end{prop}
\begin{proof}
By assumption $I'=I_N=\set{f\in\smooth{M}\,\vbar\,f(N)=\{0\}}$. Consider the ideal in $End(E)$ generated by
$I_N$. As it was mentioned above, this ideal coincides with the set of such endomorphisms
$\Phi\in End(E)$ that $\Phi|_N\equiv0$. As $I_N\cdot\mathbf{1}\subset I$, we have that
$I_N\cdot End(E)\subset I$. For any $x\in N$, consider the evaluation map
$$
\delta_x:End(E)\To End(E_x),\quad\delta_x(\Phi)=\Phi(x)
$$
As $\delta_x$ is an epimorphism, the
image of the ideal $I$ by $\delta_x$ is an ideal in $End(E_x)$. But any ideal in $End(E_x)$ is
the trivial one -- $\{0\}$, or the entire $End(E_x)$. If $\delta_x(I)=End(E_x)$, then there is an
endomorphism $\Phi\in I$, such that $\Phi_x=1$. Which implies that, there exists such neighborhood
$U$ of the point $x$ in $M$, that $\Phi_u$ is invertible for each $u\in U$. Consider two smaller
neighborhoods of $x$: $W\subset V\subset U$, and a function $f\in\smooth{M}$, such that
$f(W)=\{1\}$ and $f(M\backslash V)=\{0\}$. Construct the endomorphism
$$
\Psi=
\begin{cases}
f\cdot\Phi^{-1}\textrm{ on }V\\
0\textrm{ on }M\backslash V
\end{cases}
$$
We have that $\Psi\Phi=f\cdot\mathbf{1}$. But as $\Phi\in I$, we obtain that $f\cdot\mathbf{1}\in I$
which contradicts to the assumption that $(\smooth{M}\cdot\mathbf{1})\cap I=I_N$, because
$f(x)=1\imply f\notin I_N$. Hence, the assumption that $\delta(I)=End(E_x)$ is false. Therefore
for any $\Phi\in I$ and any $x\in N$, the value $\Phi_x$ is $0$. This, itself, implies that
$I\subset I_N\cdot End(E)$. This together with $I_N\cdot End(E)\subset I$, gives the equality
$I=I_N\cdot End(E)$. We obtain that any such ideal $I\subset End(E)$ that
$I\cap\smooth{M}\cdot\one$ is a submanifold ideal in $\smooth{M}\cdot\one$, for some compact, closed submanifold $N$, is of the form
$I=I_N\cdot End(E)$. But as we have shown, such ideals in $End(E)$ are submanifold ideals.
\end{proof}
\bigskip
\textbf{Problem: }\emph{Is any submanifold ideal in the algebra $End(E)$ of the form
$I\cdot End(E)$ where $I$ is a submanifold ideal in $\smooth{M}$?}

\begin{defn}
For any ideal $I\subset End(E)$, we call the set of points
$$
M_I=\set{x\in M\,\vbar\,\Phi_x=0,\,\all\Phi\in I}
$$
the $0$-set of the ideal $I$.
\end{defn}
\begin{prop}
If an ideal $I\subset End(E)$ is such that $M_I=\emptyset$ and the manifold $M$ is locally compact, then
$I=End(E)$.
\end{prop}
\begin{proof}
For any point $x\in M$, consider the evaluation map
$$
\delta_x:End(E)\To End(E_x),\quad\delta_x(\Phi)=\Phi_x
$$
This map is an epimorphism of associative algebras. Therefore $\delta_x(I)$ is an ideal in $End(E_x)$, which can be
only $0$ or the entire $End(E_x)$. By assumption, the $0$-set of the ideal $I$ is empty. Therefore
$\delta_x(I)=End(E_x)$. Let $\set{\Phi_x^1,\ldots,\Phi_x^m}$ be any basis of the complex vector
space $End(E_x)$, and $\set{\Phi^1,\ldots,\Phi^m}$ be such elements from $I$, that
$\delta_x(\Phi^i)=\Phi_x^i,\,i=1,\ldots,m$. As the system $\set{\Phi^1,\ldots,\Phi^m}$ is linearly
independent at the point $x$, there exists such neighborhood $U$ of $x$ in $M$, that the system
$\set{\Phi^1|_U,\ldots,\Phi^m|_U}$ is also linearly independent over $\smooth{U}$ and therefore forms a basis for
$End(E_U)$. Let $\cup U_k=M$ be such locally finite covering of $M$ that for any $U_k$ we have
a system $\set{\Phi^1_k,\ldots,\Phi^m_k}\subset I$ such that
$\set{\Phi^1_k|_{U_k},\ldots,\Phi^m_k|_{U_k}}$ is a local basis of $End(E_{U_k})$. For any
$\Psi\in End(E)$, we have: $\Psi_{U_k}=\jami{i}{}f_i^k\Phi^i_k|_{U_k}$, where
$\set{f_1^k,\ldots,f_m^k}\subset\smooth{U_k}$. Let $\set{\varphi_k}$ be a partition of unity
corresponding to the covering $\set{U_k}$. Consider the endomorphisms
$$
\widetilde{\Psi}_k=\jami{i}{}\widetilde{f}^k_i\Psi^i_k,\textrm{ where }\widetilde{f}^k_i=
\begin{cases}
\varphi_kf_i^k & \textrm{ on }U_k\\
0 & \textrm{ on }M\backslash U_k
\end{cases}
$$
As $\Phi^i_k\in I$, we have that $\widetilde{\Psi}_k\in I$, but on the other hand:
$\jami{k}{}\widetilde{\Psi}_k=\Psi$. Hence, we obtain that $\Psi\in I$.
\end{proof}

It follows from the definition of the $0$-set, that if $I_1\subseteq I_2$ then
$M_{I_2}\subseteq M_{I_1}$. Therefore for any subset $S\subset M$, the ideal $I_S\cdot End(E)$,
where
$$
I_S=\set{f\in\smooth{M}\,\vbar\,f|_S\equiv 0}
$$
is the greatest ideal with 0-set equal to $S$. From
this follows that maximal proper ideals in $End(E)$ are the ideals of the form
$I_x\cdot End(E)\cong\set{\Phi\in End(E)\,\vbar\,\Phi_x=0}$, where $x$ is a point in $M$. As it
was mentioned, such ideal is a submanifold ideal and the corresponding submanifold algeba is
$End(E)/(I_x\cdot End(E))\cong End(E_x)$ (see \cite{Masson1}).
\section{Degenerate Poisson Structures in Noncommutative Geometry}\label{sectionC}
As before, let $A$ be a unital associative complex or real algebra. A Poisson structure on $A$ is
defined as a Lie bracket $\{\,,\,\}:A\wedge A\To A$, which is a biderivation (see \cite{Dub2}, \cite{Dub3}):
$$
\brck{a}{bc}=b\brck{a}{c}+\brck{a}{b}c,\quad\all\set{a,b,c}\subset A
$$
The pair $(A,\brck{\ }{\ })$ is called a Poisson algebra. For any element $a\in A$, the derivation
$A\ni x\mapsto\brck{a}{x}\in A$ is called the Hamiltonian derivation (or just Hamiltonian)
corresponding to the element $a$ and we denote this derivation by $ham(a)$.

One method of defining a Poisson structure, which has its well-known classical analogue, is the
method which uses a symplectic form (see \cite{Dub2}, \cite{Dub3}). Let $\omega$ be a 2-form in the derivation-based differential
calculus for $A$; i.e., $\omega$ is a $Z(A)$-bilinear antisymmetric mapping:
$\omega:Der(A)\times Der(A)\To A$. Such form is said to be \emph{nondegenerate} if for any element
$a\in A$, there exists a derivation $X_a\in Der(A)$ such that $i_{X_a}\omega=-da$. The 1-forms
$i_{X_a}\omega\textrm{ and }da:Der(A)\To A$ are defined as
$(da)(X)=X(a)\textrm{ and }(i_{X_a}\omega)(X)=\omega(X_a,X),\,\all X\in Der(A)$.
It is easy to verify that if the form $\omega$ is nondegenerate then the mapping $X\mapsto i_X\omega$ from $Der(A)$
to the space of 1-forms is a monomorphism.

After this, we can state that if $\omega$ is nondegenerate then for any $a\in A$, the derivation
$X_a$ is unique. For a nondegenerate $\omega$ define an antisymmetric bracket on $A$:
$\brck{a}{b}=\omega(X_a,X_b)$. it can be verified by direct calculations that this bracket is a
biderivation and for any $a,b,c\in A$ we have (see \cite{Dub2}):
$$
\brck{\brck{a}{b}}{c}+\brck{\brck{b}{c}}{a}+\brck{\brck{c}{a}}{b}=(d\omega)(X_a,X_b,X_c)
$$
Therefore if the form $\omega$ is closed, then the Jacoby identity for this bracket is true and
so the pair $(A,\brck{\ }{\ })$ is a Poisson algebra. In this case the derivation $X_a$ is the same as
the Hamiltonian corresponding to $a$. When the $Z(A)$-module generated by the
set $\set{ham(a)\,\vbar\,a\in A}$ coincides with the entire space $Der(A)$, the two
conditions: Jacoby identity for $\brck{\ }{\ }$ and $d\omega=0$, are equivalent. A nondegenerate
and closed 2-form $\omega$ is called a \emph{symplectic form} and the pair $(A,\omega)$ is called the
corresponding \emph{symplectic structure}, or symplectic algebra. Let us formulate a
noncommutative analogue of the notion of degenerate/nondegenerate Poisson structure. The classical situation,
when $A$ is an algebra of $C^\infty$-class functions on some $C^\infty$-class manifold, dictates the following
\begin{defn}\label{def_nondeg-poisson}
A poisson bracket on an associative algebra $A$ is said to be nondegenerate if and only if the
$Z(A)$-submodule generated by the Hamiltonian derivations in $Der(A)$ coincides with the entire
$Der(A)$. Otherwise, the Poisson structure is said to be degenerate.
\end{defn}
A noncommutative generalization of a symplectic leaf for a degenerate Poisson structure can be
obtained by using of the notion of Poisson ideal.
\begin{defn}
A two-side ideal $I$ in $A$ is said to be a Poisson ideal if $\brck{I}{A}\subset I$, i.e., $I$ is
also a Lie algebra ideal. A Poisson ideal is said to be \textbf{locally proper} if
$Der(A,I)=\set{X\in Der(A)\,\vbar\,X(I)\subset I}\neq Der(A)$.
\end{defn}
If $I$ is also a submanifold ideal in $A$ (i.e., the mapping $\pi:Der(A,I)\To Der(A/I)$ is an
epimorphism) then the submanifold algebra $A/I$ can be regarded as a noncommutative analogue of a
symplectic submanifold of a Poisson manifold.

Using the above definition we formulate the following criteria for the degeneracy of a Poisson
structure on $A$.
\begin{prop}
If a Poisson algebra $A$ contains at least one locally proper Poisson ideal then the Poisson structure on $A$ is
degenerate.
\end{prop}
\begin{proof}
The $Z(A)$-module generated by $ham(A)$ is the set of the elements of the form
$\jami{i=1}{k}z_i\cdot ham(a_i)$ for some $z_i\in Z(A),\,a_i\in A,\quad 1\leq i\leq k$. It follows
from the definition of locally proper Poisson ideal, that there exists at least
one $X\in Der(A)$, such that $X(I)\nsubseteq I$. If $X\in Z(A)\cdot ham(A)$ then
$X=\sum z_i\cdot ham(a_i)$ and for any $u\in I$ we have
$X(u)=\sum z_i\cdot\underbrace{\brck{a_i}{u}}_{\in I}\in I$. This contradicts to the condition
$X(I)\nsubseteq I$. Therefore the module $Z(A)\cdot ham(A)$ does not coincide with the entire
$Der(A)$ and hence, the Poisson bracket $\brck{\ }{\ }$ is degenerate in the sense of the
definition \ref{def_nondeg-poisson}.
\end{proof}

In the classical case, when $A=\smooth{M}$, the above criteria is a necessary and sufficient condition, because if a
Poisson structure on a smooth manifold $M$ is degenerate, then we have a foliation of $M$ on
symplectic submanifolds. For any such submanifold $N$, with $\textrm{codim}>0$, the ideal of the smooth functions on $M$
vanishing on $N$ is a locally proper Poisson ideal in the Poisson algebra $\smooth{M}$.

In the classical differential calculus there is a one-to-one correspondence,
between the family of nondegenerate Poisson structures on $M$ and the family of symplectic forms
on $M$. In the noncommutative case we have the following
\begin{prop}
Let $\brck{\ }{\ }$ be a nondegenerate Poisson bracket on an associative algebra $A$ (i.e.,
$Z(A)\cdot ham(A)=Der(A)$). Then there exists a symplectic form $\omega$ such that the Poisson
bracket on $A$ defined by $\omega$ coincides with $\brck{\ }{\ }$.
\end{prop}
\begin{proof}
Let $X,Y\in Der(A)$. As the Poisson bracket $\brck{\ }{\ }$ is nondegenerate, by the definition
\ref{def_nondeg-poisson} we have $Z(A)\cdot ham(A)=Der(A)$. Therefore we have
$X=\jami{i=1}{m}u_i\cdot ham(x_i)$ and $Y=\jami{j=1}{n}v_j\cdot ham(y_j)$, for some
$\set{u_1,\ldots,u_m,v_1,\ldots,v_n}\subset Z(A)$ and
$\set{x_1,\ldots,x_m,y_1,\ldots,y_n}\subset A$. Define the value of a form $\omega$ on the pair of
vectors $(X,Y)$ as $\omega(X,Y)=\jami{i,j}{}u_iv_j\brck{x_i}{y_j}$. Let us verify that the right
side of the equality is independent of the choice of the representations of $X$ and $Y$. Let
$X=\jami{i=1}{k}u'_i\cdot ham(x'_i)$ be another representation of $X$. Then for any $1\leq j\leq n$
we have
$$
\jami{i}{}u'_iv_j\brck{x'_i}{y_j}=\jami{i}{}v_ju'_i\brck{x'_i}{y_j}=
v_jX(y_j)=\jami{p}{}v_ju_p\brck{x_p}{y_j}
$$
which implies the independence of the expression $\jami{i,j}{}u_iv_j\brck{x_i}{y_j}$ from the
representation of $X$. As the form $\omega$ is antisymmetric, from this follows the independence
from the representation of $Y$.\\
The identity of the bracket defined by the form $\omega$ and the original one is a tautological
result of the definition.\\
The form is closed because of the equality
$$
(d\omega)(ham(a),ham(b)ham(c)))=\brck{\brck{a}{b}}{c}+\brck{\brck{b}{c}}{a}+\brck{\brck{c}{a}}{b}
$$
and the assumption that $Z(A)\cdot ham(A)=Der(A)$.
\end{proof}

When $A=\smooth{M}$ for some smooth manifold $M$, the converse is also true: a Poisson
structure defined by some symplectic form is nondegenerate in the sense that
$\smooth{M}\cdot ham(M)=\fields{M}$, where $ham(M)$ denotes the space of Hamiltonian vector fields
on $M$.
\section{Poisson Structures for the Endomorphism Algebra of a Vector Bundle}\label{sectionD}
In this section we describe the family of all Poisson structures and the noncommutative analogue of the
corresponding symplectic foliations for the endomorphism algebra of a finite-dimensional vector bundle.
We start from the description of Poisson structures on the matrix algebra.

For $n\in\mathbb{N}$, let $\matr{n}$ be the algebra of $n\times n$ complex matrices. As
on any associative algebra, there is a natural Poisson structure on $\matr{n}$ defined by the
multiplicative commutator: $\brck{a}{b}=[a,b]=ab-ba$. This Poisson structure is
nondegenerate because the Hamiltonian corresponding to an element $a\in\matr{n}$ is the inner
derivation: $ham(a)=ad(a)=[a,\,\cdot]$, and all the derivations of the algebra $\matr{n}$ are
inner: $Der(\matr{n})=Int(\matr{n})$.

We have a family of Poisson structures on $\matr{n}$ parametrized by $\Complex$:
$$
\brck{a}{b}_k=k\cdot[a,b],\quad k\in\Complex
$$
All of these brackets are nondegenerate except of the case when $k=0$.
\begin{prop}
Any Poisson bracket on the algebra $\matr{n}$ is of the form
$$
\brck{a}{b}_k=k\cdot[a,b],\quad\all a,b\in\matr{n}
$$
for some $k\in\Complex$.
\end{prop}
\begin{proof}
Let $p\in\matr{n}$ be a projector: $p^2=p$. There is a one-to-one correspondence between the
subset of all projectors in $\matr{n}$ and the set of decompositions of the vector space $\Complex^n$
into a pair of complementar subspaces: $\Complex^n=X\oplus Y$. For such decomposition the
corresponding projector $p$ is the projection operator on $X$ along $Y$ and the operator $1-p$ is
the projection operator on $Y$ along $X$. As the Poisson structure is a biderivation,
for any $a\in\matr{n}$ we have the following
$$
\brck{a}{p}=\brck{a}{p^2}=p\brck{a}{p}+\brck{a}{p}p\quad\imply\quad\brck{a}{p}(1-p)=p\brck{a}{p}
$$
Consider the decomposition of $\Complex^n$, corresponding to the projector $p$:
$$
\Complex^n=\image{p}\oplus\image{1-p}
$$
If $x\in\image{p}$, then we have
$$
p\brck{a}{p}x=\brck{a}{p}\underbrace{(1-p)x}_0=0\quad\imply\quad\brck{a}{p}x\in\image{1-p}
$$
If $y\in\image{1-p}$, then
$$
\brck{a}{p}\underbrace{(1-p)y}_y=p\brck{a}{p}y\quad\imply\quad
p\brck{a}{p}y=\brck{a}{p}y\quad\imply\quad\brck{a}{p}y\in\image{p}
$$
From these we obtain that the operator $\brck{a}{p}$ maps the subspace $\image{p}$ into the
subspace $\image{1-p}$ and the subspace $\image{1-p}$ into $\image{p}$. This implies that
for any $a\in\matr{n}$, the decomposition of the operator $\brck{p}{a}$, corresponding to the
decomposition $\Complex^n=\image{p}\oplus\image{1-p}$ is of the form
$$
\brck{p}{a}=
\begin{pmatrix}
  0 & b_1 \\
  b_2 & 0
\end{pmatrix}
$$
for some $b_1:\image{1-p}\To\image{p}$ and $b_2:\image{p}\To\image{1-p}$. By definition of a
Hamiltonian derivation and because all derivations of $\matr{n}$ are inner, we have
$$
\brck{p}{a}=ham(p)a=[ham(p),a],\quad\all a\in\matr{n}
$$
Let $
ham(p)=
\begin{pmatrix}
  p_1 & q_1 \\
  p_2 & q_2
\end{pmatrix}
$
for the decomposition $\Complex^n=\image{p}\oplus\image{1-p}$. Then, for
a=
$
\begin{pmatrix}
  0 & 1 \\
  0 & 0
\end{pmatrix}
$
we obtain
$$
\brck{p}{a}=
\begin{pmatrix}
  0 & b_1 \\
  b_2 & 0
\end{pmatrix}
=
\begin{pmatrix}
  -p_2 & p_1-q_2 \\
     0 & p_2
\end{pmatrix}
\quad\imply\quad p_2=0
$$
Similarly, for
$
a=
\begin{pmatrix}
  0 & 0 \\
  1 & 0
\end{pmatrix}
$
we obtain
$$
\begin{pmatrix}
  0 & b_1 \\
  b_2 & 0
\end{pmatrix}
=
\begin{pmatrix}
  q_1     & 0 \\
  q_2-p_1 & -q_1
\end{pmatrix}
\quad\imply\quad q_1=0
$$
So, we have\
$
ham(p)=
\begin{pmatrix}
  p_1 & 0 \\
    0 & q_2
\end{pmatrix}
$.
Then, for any
$
a=
\begin{pmatrix}
  x & 0 \\
  0 & y
\end{pmatrix}
$
we obtain
$$
\brck{p}{a}=
\begin{pmatrix}
  0 & b_1 \\
  b_2 & 0
\end{pmatrix}
=
\begin{pmatrix}
  [p_1,x] & 0 \\
  0       & [q_2,y]
\end{pmatrix}
\quad\imply\quad [p_1,x]=[q_2,y]=0
$$
which implies that $p_1$ and $q_2$ are just scalars, and $ham(p)$ is of the form
$$
ham(p)=
\begin{pmatrix}
  k_1\cdot\one & 0 \\
    0 & k_2\cdot\one
\end{pmatrix},\quad k_1,k_2\in\Complex
$$
Rewrite the latter expression for $ham(p)$ in the form
$$
ham(p)=
\begin{pmatrix}
  (k_1-k_2)\cdot\one & 0 \\
        0 & 0
\end{pmatrix}
+
\begin{pmatrix}
  k_2\cdot\one & 0 \\
    0 & k_2\cdot\one
\end{pmatrix}
$$
As the matrix
$
\begin{pmatrix}
  k_2\cdot\one & 0 \\
    0 & k_2\cdot\one
\end{pmatrix}
$
is just a scalar: $k_2\cdot\mathbf{1}$, it has no effect in the commutator. Hence we can conclude
that the hamiltonian mapping on the subset of projectors is of the form
$$
ham(p)=
\begin{pmatrix}
  f(p)\cdot\one & 0 \\
     0 & 0
\end{pmatrix}
=f(p)\cdot p
$$
where $f$ is a complex-valued function on the subset of projectors in $\matr{n}$. From this easily
follows that if two projectors $p$ and $q$ does not commute then the equality
$\brck{p}{q}=[ham(p),q]=-[ham(q),p]$ implies that $f(p)=f(q)$. Now consider the set of elementary
matrices: $E_{ij}\in\matr{n},\,\set{i,j}\subset\set{1,\ldots,n}$, which forms a basis for
$\matr{n}$. Construct the following set of projectors: $\forall\set{i,j}\subset\set{1,\ldots,n}$
$$
P_{ij}=
  \begin{cases}
   E_{ij}+E_{jj}& \textrm{if }i\neq j, \\
   E_{ii}       & \textrm{if }i=j.
  \end{cases}
$$
It is clear that the set $\set{P_{ij}}$ is also a basis of the vector space $\matr{n}$. For any
$i\neq j$ we have the following commutation relations: $[P_{ij},P_{jj}]=E_{ij}\neq 0$ and $[P_{ij},P_{ii}]=-E_{ij}\neq 0$.
Therefore for any two $P_{mk}$ and $P_{kl}$, where $m\neq k$ and $k\neq l$, we have $P_{kk}$, which
does not commute with both of them. This implies that $f(P_{mk})=f(P_{kk})=f(P_{kj})$. by using a
serie of such equalities, we can obtain that
$f(P_{ij})=\textbf{const}\equiv k\in\Complex,\,\set{i,j}\subset\set{1,\ldots,n}$. As
$\set{P_{ij}}$ forms a basis of $\matr{n}$ and $ham$ is a linear map, we obtain that
$ham(T)=k\cdot T,\ \all T\in\matr{n}$.\\
Hence we can conclude that \emph{the family of Poisson structures on the algebra of matrices $\matr{n}$
is parametrized by $\Complex$ and each Poisson bracket is of the form
$\brck{S}{T}_k=k\cdot [S,T],\,k\in\Complex$}.
\end{proof}

As in the section \ref{sec_submanifolds-for-endomorphism-algebra}, let $End(E)$ be the
endomorphism algebra of a finitedimensional complex or real vector bundle $\pi:E\To M$.
As in the case of the matrix algebra, we denote by $ad(\Phi),\,\Phi\in End(E)$ the Hamiltonian derivations for
the Poisson bracket defined by the multiplicative commutator.
In this case the $\smooth{M}$-module generated by the Hamiltonian derivations is not the
entire $Der(End(E))$, because a derivation of the type $D_X,\,X\in\fields{M}$ (see the formula \ref{form_connection-expand-endom}),
defined by using of some connection on $E$, is not Hamiltonian. This follows from the fact that for any hamiltonian
derivation $ham(\Phi)=[\Phi,\,\cdot]$ and an endomorphism of $f\cdot\one,\,f\in\smooth{M}$, we
have the following
$$
ham(\Phi)(f\cdot\one)=0\textrm{ and }D_X(f\cdot\one)=X(f)\cdot\one
$$
Moreover, there are proper Poisson ideals in the Poisson algebra $(End(E),[\ ,\ ])$. Such ideal can be
constructed by a proper subset $S\subset M$, for example when $S$ consists of only one point.
\begin{lem}
For any point $x_0\in M$, the ideal $I_{x_0}(End(E))=I_{x_0}\cdot End(E)$ in the algbra $End(E)$,
generated by $I_{x_0}=\set{f\in\smooth{M}\,\vbar\,f(x_0)=0}$ is a locally proper Poisson ideal in $End(E)$.
\end{lem}
\begin{proof}
The inclusion $\brck{I_{x_0}(End(E))}{End(E)}\subset I_{x_0}(End(E))$ easily follows from the
definition of the ideal $I_{x_0}(End(E))$. Now, let us verify that there exists such derivation
$U\in Der(End(E))$ that $U(I_{x_0}(End(E)))\nsubseteq I_{x_0}(End(E))$. Let $f\in\smooth{M}$ is
such that $f(x_0)=0$ and for some tangent vector $v\in T_{x_0}(M)$: $f'(x_0)v\neq 0$. Choose such
$\Phi\in End(E)$ that $\Phi(x_0)\neq 0$. It is clear that $f\cdot\Phi\in I_{x_0}(End(E))$. Let $X$
be such vector field on $M$ that $X(x_0)=v$. For some connection on $E$ consider the derivation
$D_X:End(E)\To End(E)$. For these data we have the following
$$
D_X(f\cdot\Phi)=X(f)\cdot\Phi+f\cdot D_X(\Phi)
$$
which implies that $D_X(f\cdot\Phi)_{x_0}\neq 0$ and therefore
$D_X(f\cdot\Phi)\notin I_{x_0}(End(E))$.
\end{proof}

In the case when we have a compact, closed submanifold $N\subset M$, the ideal
$I_N(End(E))=I_N\cdot End(E)$ is a submanifold ideal (see \ref{prop_submanifold_in_base_total}).
At the same time, the quotient $End(E)/(I_N\cdot End(E))\cong End(E_N)$ is a Poisson algebra under
the bracket induced from $End(E)$. This situation can be considered as a noncommutative Poisson
submanifold in a Poisson manifold. When $N=x_0$ is a Point, the corresponding ideal
$I_{x_0}(End(E))$ is maximal and the quotient $End(E)/(I_{x_0}\cdot End(E))\cong End(E_{x_0})$ is
a nondegenerate Poisson algebra with the bracket induced from $End(E)$, and can be regarded as a
noncommutative analogue of a symplectic leaf of a Poisson structure.

Before we start the description of the family of Poisson structures on the endomorphism algebra
let us concern some general facts about Poisson algebras.

Let $(A,\brck{\ }{\ })$ be a Poisson algebra. As the center of the algebra $A$ is invariant for
any derivation, and the bracket $\brck{\ }{\ }$ is a biderivation, we have that the center of $A$
is a Lie algebra ideal for the Poisson bracket: $\brck{Z(A)}{A}\subset Z(A)$.
\begin{lem}
Let $Comm(A)$ be the minimal subalgebra of the associative algebra $A$ containing the commutator
$[A,A]=\set{[a,b]=ab-ba\,\vbar\,a,b\in A}$ (in fact
$Comm(A)$ is the set of finite sums of the elements of the type $\prod\limits_{i=1}^k[a_i,b_i]$).
For any Poisson structure on $A$, the Poisson bracket of the elements of the center of $A$ with the elements of
$Comm(A)$ is equal to $0$: $\brck{Z(A)}{Comm(A)}=\set{0}$.
\end{lem}
\begin{proof}
Because of the Leibniz rule it is sufficient to verify the statement of the lemma for the elements of the type $[a,b],\,a,b\in A$.
For any $z\in Z(A)$ we have the following
$$
\brck{z}{[a,b]}=\brck{z}{ab}-\brck{z}{ba}=a\underbrace{\brck{z}{b}}_{\in Z(A)}+
\underbrace{\brck{z}{a}}_{\in Z(A)}b-
b\underbrace{\brck{z}{a}}_{\in Z(A)}-\underbrace{\brck{z}{b}}_{\in Z(A)}a=0
$$
\end{proof}
\begin{cor}
If $Comm(A)=A$ then the multiplicative center $Z(A)$ is also in the center of the Lie algebra
corresponding to the Poisson bracket.
\end{cor}
\begin{prop}\label{prop_ham1}
If $Comm(A)=A$ and every derivation $X\in Der(A)$ such that $X(Z(A))=0$, is inner derivation, then
any Poisson structure on $A$ is defined by a linear mapping $f:A\To A/Z(A)$ via the equality
$\brck{a}{b}=[f(a),b],\quad a,b\in A$; and the mapping $f$ is $Z(A)$-linear.
\end{prop}
\begin{proof}
As it follows from the previous corollary, if $Comm(A)=A$, we have that for any $a\in A$
$$
\brck{a}{Z(a)}=\{0\}\imply ham(a)Z(A)=\{0\}\imply ham(a)\in Int(A)\cong A/Z(A)
$$
In fact, in this case the map $f:A\To A/Z(A)$ is the Hamiltonian map.

For $z\in Z(A)$ and $a\in A$, we have
$$
\begin{array}{l}
[f(za),x]=[za,f(x)]=z[a,f(x)]=z[f(a),x]=[zf(a),x]\imply\\
\\
\imply [f(za)-zf(a),x]=0,\quad\all x\in A\quad\imply f(za)-zf(a)\in Z(A)\imply\\
\\
\imply\quad [f(za)]=[zf(a)]=z[f(a)]\textrm{ in the quotient }A/Z(A)
\end{array}
$$
which implies that $f$ is $Z(A)$-linear.
\end{proof}

The antisymmetric property and the Jacoby identity for $\brck{\ }{\ }$ gives the following
properties of the mapping $f$: $\all\set{x,y,z}\subset A$
\begin{enumerate}
\item
$[f(x),y]=[x,f(y)]$
\item
$[f(x),[f(y),z]]+[f(y),[f(z),x]]+[f(z),[f(x),y]]=0$\\
or equivaletly:\\
$[[f(x),f(y)],z]+[[f(y),f(z)],x]+[[f(z),f(x)],y]=0$\\
or using the previous property:\\
$[[f^2(x),y],z]+[[f^2(y),z],x]+[[f^2(z),x],y]=0$
\end{enumerate}

Let us summarize by the following
\begin{prop}\label{prop_ham2}
Let $A$ be an associative algebra such that $Comm(A)=A$ and any derivation of $A$ which vanishes
on the center of $A$ is inner. Then there is a one-to-one correspondence between the family of
Poisson brackets on $A$ and the family of such mappings (Hamiltonians) $f:A\To A/Z(A)$ which satisfy the
conditions
\begin{itemize}
\item[\textbf{H1.}]
$f$ is a homomorphism of $Z(A)$-modules
\item[\textbf{H2.}]
$[f(a),b]=[a,f(b)],\quad\all a,b\in A$
\item[\textbf{H3.}]
$[[f^2(a),b],c]+[[f^2(b),c],a]+[[f^2(c),a],b]=0,\quad\all a,b,c\in A$
\end{itemize}
\end{prop}
It is clear that any mapping of the type $f(x)=c\cdot x$ for fixed $c\in Z(A)$ satisfies the above
conditions.

\bigskip
\bigskip
Let us apply the above result to the case when $A$ is the endomorphism algebra of a finite-
dimensional complex (or real) vector bundle $\pi:E\To M$. First, recall that any derivation of the
algebra $End(E)$, which vanishes on the center $Z(End(E))\cong\smooth{M}\cdot\one$, is inner.
Then, for the matrix algebra $\matr{n}$, when $n\geq 2$, we have that $Comm(\matr{n})=\matr{n}$. This easily
follows from the fact that any elementary matrix $E_{ij}$ is a commutator, when $i\neq j$:
$E_{ij}=[E_{ij},E_{jj}]$, and is a product of two commutators when $i=j$:
$E_{ii}=E_{ik}\cdot E_{ki}=[E_{ik},E_{kk}]\cdot [E_{ki},E_{ii}]$ for some $k\neq i$. Therefore
the same is true for the endomorphism algebra:
$$
Comm(End(E))=End(E)\textrm{, when }\dim(\pi^{-1}(x))\geq 2
$$
After these, from the Proposition \ref{prop_ham2} follows that any Poisson bracket on $End(E)$ is of the form
$\brck{\Phi}{\Psi}=[f(\Phi),\Psi]$, where
$$
f:End(E)\To End(E)/\smooth{M}
$$
is a homomorphism of $\smooth{M}$-modules. The module $End(E)/\smooth{M}$ is canonically isomorphic to the module of
sections of the subbundle of $End(E)$ consisting of the traceless endomorphisms. As $f$ is a
$\smooth{M}$-module homomorphism it is induced by by some homomorphism of vector bundles. Hence,
$f$ is pointwise: $f(\Phi)_x=f_x(\Phi_x),\,\all x\in M$. Therefore, it induces a Poisson structure
on each fiber of the endomorphism bundle:
$\brck{\Phi_x}{\Psi_x}=\brck{\Phi}{\Psi}_x=[f_x(\Phi_x),\Psi_x]$. As it was shown, any Poisson
bracket on the matrix algebra is of the form $\brck{a}{b}=k\cdot[a,b],\,k\in\Complex$, and any
Hamiltonian is of the form $f(a)=k\cdot a$. Therefore, the mapping
$$
f:End(E)\To End(E)/\smooth{M}
$$
is of the form $f(\Phi)=\lambda\cdot[\Phi]$, where $\lambda\in\smooth{M}$.

To summarize, we formulate the following \textbf{(the full description of the family of Poisson
structures on the endomorphism algebra of a vector bundle with $\dim(fiber)\geq 2$)}
\begin{prop}
Any Poisson bracket on the endomorphism algebra of a finite-dimensional complex (or real) vector
bundle, with fiber the dimension of which is $\geq 2$, is of the form
$$
\brck{\Phi}{\Psi}_\lambda=\lambda\cdot[\Phi,\Psi],\quad\all\Phi,\Psi\in End(E)
$$
where $\lambda\in\smooth{M}$.
\end{prop}
Hence, the family of Poisson structures on the endomorphism algebra $End(E)$ is parametrized by $\smooth{M}$.
\begin{cor}
Any Poisson structure on the endomorphism algebra of a finite-dimensional complex (or real) vector
bundle with $\dim(fiber)\geq 2$ is degenerate.
\end{cor}

If the fiber of a vector bundle is 1-dimensional, then the endomorphism algebra coincides with the
commutative algebra of smooth functions on the base manifold, and the Poisson structures on this algebra
are given by involutive bivector fields on the base manifold (see \cite{Lichn}).
\bibliographystyle{amsplain}

\end{document}